\documentclass{amsart}
\usepackage{amssymb}
\theoremstyle{plain}
 \newtheorem{thm}{Theorem}
 \newtheorem{cor}{Corollary}
 \newtheorem{lem}{Lemma}
 \newtheorem{prop}{Proposition}
\theoremstyle{definition}

\theoremstyle{remark}

\newcommand{\NaturalNumber}{\mathbb N}

\newcommand{\RationalNumber}{\mathbb Q}
\renewcommand{\labelenumi}{(\roman{enumi})}

\begin{document}
\title[$\bigcap_{t \geq 0}F \big( T(t) \big)
 = F \big( T(1) \big) \cap F \big( T(\sqrt2) \big)$]
{The set of common fixed points of a one-parameter continuous semigroup
 of mappings
is $F \big( T(1) \big) \cap F \big( T(\sqrt2) \big)$}
\author[T. Suzuki]{Tomonari Suzuki}
\date{}
\address{
Department of Mathematics,
Kyushu Institute of Technology,
1-1, Sensuicho, To\-bata\-ku, Kitakyushu 804-8550, Japan}
\email{suzuki-t@mns.kyutech.ac.jp}

\begin{abstract}
In this paper, we prove the following theorem:
Let $\{ T(t) : t \geq 0 \}$ be a one-parameter continuous semigroup of mappings
 on a subset $C$ of a Banach space $E$.
The set of fixed points of $T(t)$ is denoted by $F \big( T(t) \big)$
 for each $t \geq 0$.
Then
 $$ \bigcap_{t \geq 0} F \big( T(t) \big)
 = F \big( T(1) \big) \cap F \big( T(\sqrt 2) \big) $$
 holds.
Using this theorem, we discuss convergence theorems
 to a common fixed point of $\{ T(t) : t \geq 0 \}$.
\end{abstract}
\maketitle

\section{Introduction}
\label{SC:intro}

Let $C$ be a subset of a Banach space $E$, and
 let $T$ be a nonexpansive mapping on $C$, i.e.,
 $\| Tx - Ty \| \leq \| x - y \|$
 for all $x, y \in C$.
We know that $T$ has a fixed point
 in the case that $E$ is uniformly convex and
 $C$ is bounded, closed and convex;
 see Browder \cite{REF:Browder1965_ProcNAS_3},
 G\"ohde \cite{REF:Gohde1965}, and Kirk \cite{REF:Kirk1965_AMMonth}.
We denote by $F(T)$ the set of fixed points of $T$.

Let $\tau$ be a Hausdorff topology on $E$.
A family of mappings $\{ T(t) : t \geq 0 \}$ is called
 a one-parameter $\tau$-continuous semigroup of mappings on $C$
 if the following are satisfied:
 \begin{enumerate}
 \renewcommand{\labelenumi}{(sg \arabic{enumi})}
 \item
 $ T(s+t) = T(s) \circ T(t) $ for all $s, t \geq 0$;
 \item
 for each $x \in X$,
 the mapping $t \mapsto T(t)x $ from $[0,\infty)$ into $C$ is continuous
 with respect to $\tau$.
 \end{enumerate}
As topology $\tau$,
 we usually consider the strong topology of $E$.
Also, a family of mappings $\{ T(t) : t \geq 0 \}$ is called
 a one-parameter $\tau$-continuous semigroup
 of nonexpansive mappings on $C$
 (in short, nonexpansive semigroup)
 if (sg~1), (sg~2) and the following (sg~3) are satisfied:
\begin{enumerate}
\setcounter{enumi}{2}
 \renewcommand{\labelenumi}{(sg \arabic{enumi})}
\item
 for each $t \geq 0$,
 $T(t)$ is a nonexpansive mapping on $C$.
\end{enumerate}
We know that nonexpansive semigroup $\{ T(t) : t \geq 0 \}$
 has a common fixed point
 in the case that $E$ is uniformly convex and
 $C$ is bounded, closed and convex;
 see Browder \cite{REF:Browder1965_ProcNAS_3}.
Moreover, in 1974, Bruck \cite{REF:Bruck1974_Pacific} proved that
 nonexpansive semigroup $\{ T(t) : t \geq 0 \}$ has a common fixed point
 in the case that $C$ is weakly compact, convex, and
 has the fixed point property for nonexpansive mappings.

In this paper, we prove the following theorem:
Let $\{ T(t) : t \geq 0 \}$ be a one-parameter $\tau$-continuous semigroup
 of mappings
 on a subset $C$ of a Banach space $E$
 for some Hausdorff topology $\tau$ on $E$.
Then
 $$ \bigcap_{t \geq 0} F \big( T(t) \big)
 = F \big( T(1) \big) \cap F \big( T(\sqrt 2) \big) $$
 holds.
Using this theorem, we discuss convergence theorems
 to a common fixed point of nonexpansive semigroups $\{ T(t) : t \geq 0 \}$.

\section{Preliminaries}
\label{SC:preliminaries}

Throughout this paper
 we denote
 by $\RationalNumber$ the set of rational numbers, and
 by $\NaturalNumber$ the set of positive integers.
For real number $t$,
 we denote by $[t]$
 the maximum integer not exceeding $t$.
It is obvious that
 for each real number $t$,
 there exists $\varepsilon \in [0,1)$ such that
 $t = [t] + \varepsilon$.

We recall that
 a Banach space $E$ is called strictly convex if
 $\| x+y \| / 2 < 1$
 for all $ x, y \in E $ with $ \| x \| = \| y \| = 1 $
 and $ x \neq y $.
A Banach space $E$ is called uniformly convex if
 for each $\varepsilon > 0$,
 there exists $\delta > 0$ such that
 $\| x+y \| / 2 < 1 - \delta$
 for all $ x, y \in E $ with $ \| x \| = \| y \| = 1 $
 and $\| x - y \| \geq \varepsilon$.
It is clear that a uniformly convex Banach space is strictly convex.
The norm of $E$ is called Fr\'echet differentiable if
 for each $x \in E$ with $\| x \| = 1$,
 $\lim_{t \rightarrow 0}
 (\| x + t y \| - \| x \|)/t$ exists and is attained
 uniformly in $y \in E$ with $\| y \| = 1$.

The following Lemma is the corollary of
 Bruck's result in \cite{REF:Bruck1973_TransAMS}.

\begin{lem}[Bruck \cite{REF:Bruck1973_TransAMS}]
\label{LEM:Bruck}
Let $C$ be a subset of a strictly convex Banach space $E$.
Let $S$ and $T$ be nonexpansive mappings from $C$ into $E$
 with common fixed point.
Then for each $\lambda \in (0,1)$,
 a mapping $U$ from $C$ into $E$ defined by
 $ U x = \lambda S x + (1-\lambda) T x $
 for $x \in C$ is nonexpansive and
 $ F(U) = F(S) \cap F(T)$ holds.
\end{lem}

\begin{proof}
It is obvious that $F(U) \supset F(S) \cap F(T)$.
Fix $x \in F(U)$ and $w \in F(S) \cap F(T)$.
Then we have
 \begin{align*}
 \| x - w \|
 &= \| \lambda Sx + (1-\lambda) Tx - w \| \\*
 &\leq \lambda \| Sx - w \| + (1-\lambda) \| Tx- w \| \\
 &\leq \lambda \| x - w \| + (1-\lambda) \| x- w \| \\*
 &= \| x - w \|
 \end{align*}
 and hence
 $$ \| x - w \| = \| \lambda Sx + (1-\lambda) Tx - w \|
 = \| Sx - w \| = \| Tx - w \| . $$
So, from the strict convexity of $E$,
 we obtain
 $$ \lambda Sx + (1-\lambda) Tx = Sx = Tx .$$
Hence $x \in F(S) \cap F(T)$.
This completes the proof.
\end{proof}

The following four convergence theorems for nonexpansive mappings
 are well-known.

\begin{thm}[Baillon \cite{REF:Baillon1975}]
\label{THM:Baillon}
Let $C$ be a bounded closed convex subset of a Hilbert space $E$.
Let $T$ be a nonexpansive mapping on $C$.
Let $x \in C$ and
 define a sequence $\{ x_n \}$ in $C$ by
 $$ x_{n}
 = \frac{T x + T^2 x + T^3 x + \cdots + T^n x}{n} $$
 for $n \in \NaturalNumber$.
Then
 $\{ x_n \}$ converges weakly to a fixed point of $T$.
\end{thm}

\begin{thm}[Reich \cite{REF:Reich1979_JMAA}]
\label{THM:Reich}
Let $E$ be a uniformly convex Banach space
 whose norm is Fr\'echet differentiable.
Let $T$ be a nonexpansive mapping
 on a bounded closed convex subset $C$ of $E$.
Define a sequence $\{ x_n \}$ in $C$ by
 $x_1 \in C$ and
 $ x_{n+1} = \alpha_n \; Tx_n + (1-\alpha_n)\; x_n $
 for $n \in \NaturalNumber$.
 where $\{ \alpha_n \}$ is a sequence in $[0,1]$
 satisfying $\sum_{n=1}^\infty \alpha_n (1-\alpha_n) = \infty$.
Then
 $\{ x_n \}$ converges weakly to a fixed point of $T$.
\end{thm}

\begin{thm}[Browder \cite{REF:Browder1967_ARMA}]
\label{THM:Browder}
Let $C$ be a bounded closed convex subset of a Hilbert space $E$, and
 let $T$ be a nonexpansive mapping on $C$.
Let $\{ \lambda_n \}$ be a sequence in $(0,1)$ converging to $0$.
Fix $u \in C$ and
 define a sequence $\{ x_n \}$ in $C$ by
 $ x_n = (1-\lambda_n) \; T x_n + \lambda_n \; u $
 for $n \in \NaturalNumber$.
Then
 $\{ x_n \}$ converges strongly to a fixed point of $T$.
\end{thm}

\begin{thm}[Wittmann \cite{REF:Wittmann1992_ArchMath}]
\label{THM:Wittmann}
Let $C$ be a bounded closed convex subset of a Hilbert space $E$, and
 let $T$ be a nonexpansive mapping on $C$.
Let $u \in C$ and
 define a sequence $\{ x_n \}$ in $C$ by $x_1 \in C$ and
 $ x_{n+1}
 = (1-\lambda_n) \; T x_n + \lambda_n \; u $
 for $n \in \NaturalNumber$,
 where $\{ \lambda_n \}$ is a sequence in $[0,1]$ satisfying the following:
 $$ \lim_{n \rightarrow \infty} \lambda_n = 0; \quad
 \sum_{n=1}^\infty \lambda_n = \infty; \quad\text{and}\quad
 \sum_{n=1}^\infty | \; \lambda_{n+1} - \lambda_n \; | < \infty . $$
Then
 $\{ x_n \}$ converges strongly to a fixed point of $T$.
\end{thm}

\section{Lemmas}
\label{SC:lemmas}

In this Section, we prove two lemmas,
 which are used in Section \ref{SC:result}.

\begin{lem}
\label{LEM:sum}
Let $t$ be a nonnegative real number
 and let $\{ \beta_n \}$ be a sequence in $(0,\infty)$ converging to $0$.
Define sequences $\{ \delta_n \}$ in $[0,\infty)$ and
 $\{ k_n \}$ in $\NaturalNumber \cup \{ 0 \}$ as follows:
\begin{itemize}
\item $\delta_1 = t$;
\item $k_n = [\delta_n / \beta_n]$ for $n \in \NaturalNumber$;
\item $\delta_{n+1} = \delta_n - k_n \beta_n$
 for $n \in \NaturalNumber$.
\end{itemize}
Then the following hold:
\begin{enumerate}
\item $0 \leq \delta_{n+1} < \beta_n$ for all $n \in \NaturalNumber$;
\item $k_n \in \NaturalNumber \cup \{ 0 \}$ for all $n \in \NaturalNumber$;
\item $\{ \delta_n \}$ converges to $0$;
\item $\sum_{j=1}^{n} k_j \; \beta_j + \delta_{n+1} = t$
 for all $n \in \NaturalNumber$; and
\item $\sum_{j=1}^{\infty} k_j \; \beta_j = t$.
\end{enumerate}
\end{lem}

\begin{proof}
We put $\varepsilon_n \in [0,1)$ with
 $$ \frac{\delta_n}{\beta_n} = k_n + \varepsilon_n $$
 for $n \in \NaturalNumber$.
We have
 $$ \delta_{n+1} = \delta_n - k_n \beta_n
 = \varepsilon_n \; \beta_n
 < \beta_n $$
 for all $n \in \NaturalNumber$.
{}From this, we also have $\delta_{n+1} = \varepsilon_n \; \beta_n \geq 0$.
This implies (i).
It is obvious that (ii) and (iii) follow from (i).
Let us prove (iv).
We have
 $$ k_1 \; \beta_1 + \delta_2
 = k_1 \; \beta_1 + (\delta_1 - k_1 \; \beta_1)
 = \delta_1 = t .$$
We assume (iv) holds for some $n \in \NaturalNumber$.
Then we have
 \begin{align*}
 \sum_{j=1}^{n+1} k_j \; \beta_j + \delta_{n+2}
 &= \sum_{j=1}^{n+1} k_j \; \beta_j + (\delta_{n+1} - k_{n+1} \; \beta_{n+1})
 \\*
 &= \sum_{j=1}^{n} k_j \; \beta_j + \delta_{n+1} \\*
 &= t .
 \end{align*}
So, by induction, we have (iv).
{}From (iii) and (iv), we have
 $$ \sum_{n=1}^\infty k_n \beta_n = t . $$
This completes the proof.
\end{proof}

\begin{lem}
\label{LEM:Q}
Let $\alpha$ and $\beta$ be positive real numbers
 satisfying $\alpha/\beta \notin \RationalNumber$.
Define sequences $\{ \alpha_n \}$ in $(0, \infty)$ and
 $\{ k_n \}$ in $\NaturalNumber$ as follows:
\begin{itemize}
\item $\alpha_1 = \max\{ \alpha, \beta \}$;
\item $\alpha_2 = \min\{ \alpha, \beta \}$;
\item $k_n = [ \; \alpha_n / \alpha_{n+1} \; ]$
 for all $n \in \NaturalNumber$;
\item $\alpha_{n+2} = \alpha_n -  k_n \; \alpha_{n+1}$
 for all $n \in \NaturalNumber$.
\end{itemize}
Then the following hold:
\begin{enumerate}
\item $0 < \alpha_{n+1} < \alpha_n$ for all $n \in \NaturalNumber$;
\item $k_n \in \NaturalNumber$ for all $n \in \NaturalNumber$;
\item $\alpha_n / \alpha_{n+1} \notin \RationalNumber$
 for all $n \in \NaturalNumber$; and
\item $\{ \alpha_n \}$ converges to $0$.
\end{enumerate}
\end{lem}

\begin{proof}
We note that (i) implies (ii).
By the assumption of $\alpha / \beta \notin \RationalNumber$,
 we have $\alpha \neq \beta$.
Hence
 $$ \alpha_1 = \max\{ \alpha, \beta \}
 > \min\{ \alpha, \beta \} = \alpha_2 > 0. $$
It is obvious that $\alpha_1 / \alpha_2 \notin \RationalNumber$.
We assume that $0 < \alpha_{j+1} < \alpha_{j}$ and
 $\alpha_j / \alpha_{j+1} \notin \RationalNumber$
 for some $j \in \NaturalNumber$.
Since
 $\alpha_{j+2} = \alpha_j -  k_j \; \alpha_{j+1}$,
 we have
 $$ \frac{\alpha_{j+2}}{\alpha_{j+1}}
 = \frac{\alpha_j}{\alpha_{j+1}} -  k_j \notin \RationalNumber $$
 and hence $\alpha_{j+1} / \alpha_{j+2} \notin \RationalNumber$.
Put $\varepsilon_j \in [0,1)$ satisfying
 $$ \frac{\alpha_j}{\alpha_{j+1}} = k_j + \varepsilon_j .$$
Since $\alpha_j / \alpha_{j+1} \notin \RationalNumber$,
 we note that $\varepsilon_j > 0$.
We have
 $$ \alpha_{j+2} = \alpha_j -  k_j \; \alpha_{j+1}
 = \varepsilon_j \; \alpha_{j+1} < \alpha_{j+1} .$$
{}From this, we also have
 $ \alpha_{j+2} = \varepsilon_j \; \alpha_{j+1} > 0 $.
Therefore we have shown that
 $0 < \alpha_{j+2} < \alpha_{j+1}$ and
 $\alpha_{j+1} / \alpha_{j+2} \notin \RationalNumber$.
By induction, we have (i), (ii) and (iii).
Let us prove (iv).
Since $\{ \alpha_n \}$ is a sequence of positive real numbers
 and strictly decreasing,
 $\{ \alpha_n \}$ converges to some $\alpha_\infty \in [0,\infty)$.
We assume $\alpha_\infty > 0$.
Then we can choose $j \in \NaturalNumber$ such that
 $$ \alpha_\infty < \alpha_{j+1} < \alpha_j < 2 \alpha_\infty .$$
We have
 $$ k_j = \left[ \frac{\alpha_j}{\alpha_{j+1}} \right] = 1
 \quad\text{and}\quad
 \alpha_{j+2} = \alpha_j - k_j \; \alpha_{j+1}
 = \alpha_j - \alpha_{j+1} < \alpha_\infty . $$
This is a contradiction.
Therefore $\alpha_\infty = 0$ and this implies (iv).
This completes the proof.
\end{proof}

\section{Main Results}
\label{SC:result}

In this Section, we give our main results.
We know the following.

\begin{prop}
\label{PROP:cluster}
Let $E$ be a Banach space and let $\tau$ be a Hausdorff topology on $E$.
Let $\{ T(t) : t \geq 0 \}$ be
 a one-parameter $\tau$-continuous semigroup of mappings
 on a subset $C$ of $E$.
Let $\{ \alpha_n \}$ be a sequence in $[0, \infty)$
 converging to $\alpha_\infty \in [0,\infty)$,
 and satisfying $\alpha_n \neq \alpha_\infty$ for all $n \in \NaturalNumber$.
Suppose that
 $z \in C$ satisfies
 $$ T(\alpha_n) z = z $$
 for all $n \in \NaturalNumber$.
Then $z$ is a common fixed point of $\{ T(t) : t \geq 0 \}$.
\end{prop}

\begin{proof}
We note that
 $$ T(\alpha_\infty) z
 = \mathop{\tau\text{-}\lim}_{n \rightarrow \infty} T(\alpha_n)z = z . $$
We put
 $$ \beta_n = | \alpha_n - \alpha_\infty | > 0 $$
 for $n \in \NaturalNumber$.
By the assumption,
 $\{ \beta_n \}$ is a sequence in $(0,\infty)$
 converging to $0$.
Since
 $$ \max \{ \alpha_n, \alpha_\infty \}
 = \min \{ \alpha_n, \alpha_\infty \} + \beta_n , $$
 we have
 \begin{align*}
 T(\beta_n) z
 &= T(\beta_n) \circ T \big( \min \{ \alpha_n, \alpha_\infty \} \big) z \\*
 &= T \big( \beta_n + \min \{ \alpha_n, \alpha_\infty \} \big) z
 = T \big( \max \{ \alpha_n, \alpha_\infty \} \big) z \\*
 &= z
 \end{align*}
 for all $n \in \NaturalNumber$.
We also have
 $$ T(0) z = T(0) \circ T(\alpha_1) z = T(0+\alpha_1)z = T(\alpha_1) z = z . $$
Fix $t > 0$.
Then by Lemma \ref{LEM:sum},
 there exists a sequence $\{ k_n \}$ in $\NaturalNumber \cup \{ 0 \}$
 such that
 $$ \sum_{n=1}^\infty k_n \beta_n = t . $$
For each $n \in \NaturalNumber$
 with $\sum_{j=1}^n k_j \beta_j > 0$,
 we obtain
 \begin{align*}
 T \left( \sum_{j=1}^n k_j \beta_j \right) z
 &= T(\beta_n)^{k_n} \circ T(\beta_{n-1})^{k_{n-1}} \circ \cdots
  \circ T(\beta_2)^{k_2} \circ T(\beta_1)^{k_1} z \\*
 &= T(\beta_n)^{k_n} \circ T(\beta_{n-1})^{k_{n-1}} \circ \cdots
  \circ T(\beta_2)^{k_2} z \\
 &= \cdots = T(\beta_n)^{k_n} z \\*
 &= z ,
 \end{align*}
 where $T(\beta_j)^0$ is the identity mapping on $C$.
Hence, we have
 $$  T(t) z
 = \mathop{\tau\text{-}\lim}_{n \rightarrow \infty} \;
  T \left( \sum_{j=1}^n k_j \beta_j \right) z
 = z . $$
This completes the proof.
\end{proof}

We now prove one of our main results.

\begin{prop}
\label{PROP:main}
Let $E$ be a Banach space and let $\tau$ be a Hausdorff topology on $E$.
Let $\{ T(t) : t \geq 0 \}$ be
 a one-parameter $\tau$-continuous semigroup of mappings
 on a subset $C$ of $E$.
Let $\alpha$ and $\beta$ be positive real numbers
 satisfying $\alpha/\beta \notin \RationalNumber$.
Then
 $$ \bigcap_{t \geq 0} F \big( T(t) \big)
 = F \big( T(\alpha) \big) \cap F \big( T(\beta) \big) $$
 holds.
\end{prop}

\begin{proof}
It is obvious that
 $$ \bigcap_{t \geq 0} F \big( T(t) \big)
 \subset F \big( T(\alpha) \big) \cap F \big( T(\beta) \big) . $$
We fix $z \in F \big( T(\alpha) \big) \cap F \big( T(\beta) \big)$.
Define sequences $\{ \alpha_n \}$ in $(0, \infty)$ and
 $\{ k_n \}$ in $\NaturalNumber$ as in Lemma \ref{LEM:Q}.
By the assumption, we have
 $$ T(\alpha_1)z = T(\max\{ \alpha, \beta \})z = z
 \quad\text{and}\quad
 T(\alpha_2)z = T(\min\{ \alpha, \beta \})z = z . $$
If $T(\alpha_j)z=T(\alpha_{j+1})z=z$, then we have
 $$ T(\alpha_{j+2}) z
 = T(\alpha_{j+2}) \circ T(\alpha_{j+1})^{k_j} z
 = T(\alpha_{j+2} + k_j \; \alpha_{j+1}) z
 = T(\alpha_j) z
 = z .$$
So, by induction, we have
 $T(\alpha_n)z=z$ for all $n \in \NaturalNumber$.
Since $\{ \alpha_n \}$ is a positive real sequence converging to $0$,
 we have $z$ is a common fixed point of $\{ T(t) : t \geq 0 \}$
 by Proposition \ref{PROP:cluster}.
This completes the proof.
\end{proof}

As a direct consequence of Proposition \ref{PROP:main},
 we obtain the following.

\begin{cor}
\label{COR:main}
Let $E$ be a Banach space and let $\tau$ be a Hausdorff topology on $E$.
Let $\{ T(t) : t \geq 0 \}$ be
 a one-parameter $\tau$-continuous semigroup of mappings
 on a subset $C$ of $E$.
Then
 $$ \bigcap_{t \geq 0} F \big( T(t) \big)
 = F \big( T(1) \big) \cap F \big( T(\sqrt 2) \big) $$
 holds.
\end{cor}

Using Lemma \ref{LEM:Bruck},
 we obtain the following.

\begin{cor}
\label{COR:sc}
Let $E$ be a strictly convex Banach space and
 let $\tau$ be a Hausdorff topology on $E$.
Let $\{ T(t) : t \geq 0 \}$ be
 a one-parameter $\tau$-continuous semigroup of nonexpansive mappings
 on a subset $C$ of $E$.
Let $\alpha$ and $\beta$ be positive real numbers
 satisfying $\alpha/\beta \notin \RationalNumber$, and
 $F \big( T(\alpha) \big) \cap F \big( T(\beta) \big) \neq \varnothing$.
Then
 $$ \bigcap_{t \geq 0} F \big( T(t) \big)
 = \{ z \in C : \lambda T(\alpha)z + (1-\lambda) T(\beta)z = z \} $$
 holds for every $\lambda \in (0,1)$.
\end{cor}

\begin{cor}
\label{COR:uc}
Let $E$ be a uniformly convex Banach space and
 let $\tau$ be a Hausdorff topology on $E$.
Let $\{ T(t) : t \geq 0 \}$ be
 a one-parameter $\tau$-continuous semigroup of nonexpansive mappings
 on a bounded closed convex subset $C$ of $E$.
Let $\alpha$ and $\beta$ be positive real numbers
 satisfying $\alpha/\beta \notin \RationalNumber$.
Then
 $$ \bigcap_{t \geq 0} F \big( T(t) \big)
 = \{ z \in C : \lambda T(\alpha)z + (1-\lambda) T(\beta)z = z \} $$
 holds for every $\lambda \in (0,1)$.
\end{cor}

\section{Convergence Theorems}
\label{SC:convergence-theorems}

Several authors have studied about convergence theorems
 for one-parameter nonexpansive semigroups;
 see
 \cite{REF:Atsushiba_Takahashi2002_1,
 REF:Baillon1976,
 REF:Hirano1982_JMSJapan,
 REF:Miyadera_Kobayasi1982_NATMA,
 REF:Shioji_Takahashi1998_NATMA,
 REF:TS2003_ProcAMS,
 REF:Suzuki-Takahashi}
 and others.
For example, Suzuki and Takahashi prove in \cite{REF:Suzuki-Takahashi}
 the following:
Let $C$ be a compact convex subset of a Banach space $E$ and
 let $\{ T(t) : t \geq 0 \}$ be
 a one-parameter strongly continuous semigroup of nonexpansive mappings on $C$.
Let $x_1 \in C$ and
 define a sequence $\{ x_n \}$ in $C$ by
 $$ x_{n+1} =  \frac{\lambda}{t_n} \int_{0}^{t_n} T(s) x_n \; d s
  +(1 - \lambda) x_n $$
 for $n \in \NaturalNumber$,
 where $\lambda$ is a constant in $(0,1)$, and
 $ \{ t_n \}$ is a sequence in $(0,\infty)$
 satisfying $$ \lim_{n \rightarrow \infty} t_n = \infty \quad\text{and}\quad
 \lim_{n \rightarrow \infty} \frac{t_{n+1}}{t_n} = 1 . $$
Then $\{ x_n \}$ converges strongly to
 a common fixed point $z_0$ of $\{ T(t) : t \geq 0 \}$.

Using Proposition \ref{PROP:main},
 we can prove many convergence theorems to a common fixed point
 of a one-parameter continuous semigroup of nonexpansive mappings.
In this Section,
 we state some of them.
We discuss five types of convergence theorems.
Five types are the types of
 Baillon \cite{REF:Baillon1975},
 Krasnoselskii-Mann
 \cite{REF:KrasnoselskiiMA1955_Uspekhi, REF:Mann1953_ProcAMS},
 Ishikawa \cite{REF:Ishikawa1979_Pacific},
 Browder \cite{REF:Browder1967_ARMA}, and
 Halpern \cite{REF:Halpern1967_BullAMS}.
We first state the following,
 which are connected with Baillon's type iteration;
 see pages 63 and 83 in \cite{REF:Takahashi_ybook}.

\begin{thm}
\label{THM:new-Baillon}
Let $E$ be a Hilbert space and
 let $\tau$ be a Hausdorff topology on $E$.
Let $\{ T(t) : t \geq 0 \}$ be a $\tau$-continuous semigroup
 of nonexpansive mappings on a bounded closed convex subset $C$ of $E$.
Fix $\alpha, \beta > 0$ with $\alpha / \beta \notin \RationalNumber$.
Let $x \in C$ and
 define a sequence $\{ x_n \}$ in $C$ by
 $$ x_{n}
 = \frac{\displaystyle\sum_{k=1}^n
 \displaystyle\sum_{\ell=1}^n T(k \; \alpha + \ell \; \beta) \; x}{n^2} $$
 for $n \in \NaturalNumber$.
Then
 $\{ x_n \}$ converges weakly to a common fixed point of
 $\{ T(t) : t \geq 0 \}$.
\end{thm}

\begin{proof}
We note that
 $$ \sum_{k=1}^n \sum_{\ell=1}^n T(k \; \alpha + \ell \; \beta) \; x
 = \sum_{k=1}^n \sum_{\ell=1}^n
  T(\alpha)^k \circ T(\beta)^\ell \; x . $$
So, $\{ x_n \}$ converges weakly to a common fixed point $z$ of
 $T(\alpha)$ and $T(\beta)$.
Such $z$ is a common fixed point of $\{ T(t) : t \geq 0 \}$
 by Proposition \ref{PROP:main}.
This completes the proof.
\end{proof}

\begin{thm}
\label{THM:new-Baillon2}
Let $E$ be a Hilbert space and
 let $\tau$ be a Hausdorff topology on $E$.
Let $\{ T(t) : t \geq 0 \}$ be a $\tau$-continuous semigroup
 of nonexpansive mappings on a bounded closed convex subset $C$ of $E$.
Fix $\alpha, \beta > 0$ with $\alpha / \beta \notin \RationalNumber$.
Let $x \in C$ and
 define a sequence $\{ x_n \}$ in $C$ by
 $$ x_{n}
 = \frac{\displaystyle\sum_{k=1}^n
  \left( \frac{T(\alpha)+T(\beta)}{2} \right)^k x}{n} $$
 for $n \in \NaturalNumber$.
Then
 $\{ x_n \}$ converges weakly to a common fixed point of
 $\{ T(t) : t \geq 0 \}$.
\end{thm}

\begin{proof}
By Theorem \ref{THM:Baillon},
 $\{ x_n \}$ converges weakly to $z$,
 which is a fixed point of $\big( T(\alpha) + T(\beta) \big)/2$.
So, by Corollary \ref{COR:uc},
 $z$ is a common fixed point of $\{ T(t) : t \geq 0 \}$.
This completes the proof.
\end{proof}

We next state the following,
 which are connected with Krasnoselskii-Mann's type iteration;
 see Reich \cite{REF:Reich1979_JMAA} and
 Suzuki \cite{REF:TS2002_JNCA,REF:TSP_mann2uo_1_04}.

\begin{thm}
\label{THM:new-Reich}
Let $E$ be a uniformly convex Banach space
 whose norm is Fr\'echet differentiable and
 let $\tau$ be a Hausdorff topology on $E$.
Let $\{ T(t) : t \geq 0 \}$ be a $\tau$-continuous semigroup
 of nonexpansive mappings on a bounded closed convex subset $C$ of $E$.
Fix $\alpha, \beta > 0$ with $\alpha / \beta \notin \RationalNumber$,
 and $\kappa, \lambda > 0$ with $\kappa + \lambda < 1$.
Define a sequence $\{ x_n \}$ in $C$ by
 $x_1 \in C$ and
 $$ x_{n+1} = \kappa \; T(\alpha) x_n + \lambda \; T(\beta) x_n
 + (1-\kappa-\lambda) x_n , $$
 for $n \in \NaturalNumber$.
Then
 $\{ x_n \}$ converges weakly to a common fixed point of
 $\{ T(t) : t \geq 0 \}$.
\end{thm}

\begin{proof}
By Theorem \ref{THM:Reich},
 $\{ x_n \}$ converges weakly to $z$,
 which is a fixed point of
 $$ \frac{\kappa}{\kappa+\lambda} T(\alpha)
 + \frac{\lambda}{\kappa+\lambda} T(\beta) . $$
So, by Corollary \ref{COR:uc},
 $z$ is a common fixed point of $\{ T(t) : t \geq 0 \}$.
This completes the proof.
\end{proof}

\begin{thm}
\label{THM:new-Suzuki}
Let $E$ be a Banach space and
 let $\tau$ be a Hausdorff topology on $E$.
Let $\{ T(t) : t \geq 0 \}$ be a $\tau$-continuous semigroup
 of nonexpansive mappings on a compact convex subset $C$ of $E$.
Fix $\alpha, \beta > 0$ with $\alpha / \beta \notin \RationalNumber$,
 and $\lambda \in (0,1)$.
Define a sequence $\{ x_n \}$ in $C$ by
 $x_1 \in C$ and
 $$ x_{n+1} = \lambda \; \frac{\displaystyle\sum_{k=1}^n
 \displaystyle\sum_{\ell=1}^n T(k \; \alpha + \ell \; \beta) x_n}{n^2}
 + (1-\lambda) x_n , $$
 for $n \in \NaturalNumber$.
Then
 $\{ x_n \}$ converges strongly to a common fixed point of
 $\{ T(t) : t \geq 0 \}$.
\end{thm}

Using Ishikawa's result in \cite{REF:Ishikawa1979_Pacific},
 we obtain the following.

\begin{thm}
\label{THM:new-Ishikawa}
Let $E$ be a Banach space and
 let $\tau$ be a Hausdorff topology on $E$.
Let $\{ T(t) : t \geq 0 \}$ be a $\tau$-continuous semigroup
 of nonexpansive mappings on a compact convex subset $C$ of $E$.
Fix $\alpha, \beta > 0$ with $\alpha / \beta \notin \RationalNumber$
 and $\kappa, \lambda \in (0,1)$.
Define a sequence $\{ x_n \}$ in $C$ by
 $x_1 \in C$ and
 $$ x_{n+1} = \big( \lambda \; T(\alpha) + (1-\lambda) \; I \big)
 \circ \big( \kappa \; T(\beta) + (1-\kappa) \; I \big)^n x_n $$
 for $n \in \NaturalNumber$,
 where $I$ is the identity mapping on $C$.
Then
 $\{ x_n \}$ converges strongly to a common fixed point of
 $\{ T(t) : t \geq 0 \}$.
\end{thm}

We next state the following,
 which is connected with Browder's type implicit iteration.
We note that
 $$ x \mapsto (1-\lambda) \; Tx + \lambda \; u $$
 is a contractive mapping
 if $T$ is a nonexpansive mapping and $\lambda \in (0,1)$.
By the Banach contraction principle \cite{REF:Banach1922_FundMath},
 such mappings have a unique fixed point.

\begin{thm}
\label{THM:new-Browder}
Let $E$ be a Hilbert space and
 let $\tau$ be a Hausdorff topology on $E$.
Let $\{ T(t) : t \geq 0 \}$ be a $\tau$-continuous semigroup
 of nonexpansive mappings on a bounded closed convex subset $C$ of $E$.
Fix $\alpha, \beta > 0$ with $\alpha / \beta \notin \RationalNumber$.
Let $u \in C$ and
 define a sequence $\{ x_n \}$ in $C$ by
 $$ x_n
 = \frac{1-\lambda_n}{2} T(\alpha) x_n + \frac{1-\lambda_n}{2} T(\beta) x_n
 + \lambda_n \; u $$
 for $n \in \NaturalNumber$,
 where $\{ \lambda_n \}$ is a sequence in $(0,1)$ converging to $0$.
Then
 $\{ x_n \}$ converges strongly to a common fixed point of
 $\{ T(t) : t \geq 0 \}$.
\end{thm}

We finally state the following,
 which is connected with Halpern's type explicit iteration;
 see Wittmann \cite{REF:Wittmann1992_ArchMath}.

\begin{thm}
\label{THM:new-Wittmann}
Let $E$ be a Hilbert space and
 let $\tau$ be a Hausdorff topology on $E$.
Let $\{ T(t) : t \geq 0 \}$ be a $\tau$-continuous semigroup
 of nonexpansive mappings on a bounded closed convex subset $C$ of $E$.
Fix $\alpha, \beta > 0$ with $\alpha / \beta \notin \RationalNumber$.
Let $u \in C$ and
 define a sequence $\{ x_n \}$ in $C$ by
 $x_1 \in C$ and
 $$ x_{n+1}
 = \frac{1-\lambda_n}{2} T(\alpha) x_n + \frac{1-\lambda_n}{2} T(\beta) x_n
 + \lambda_n \; u $$
 for $n \in \NaturalNumber$,
 where $\{ \lambda_n \}$ is a sequence in $[0,1]$ satisfying the following:
 $$ \lim_{n \rightarrow \infty} \lambda_n = 0; \quad
 \sum_{n=1}^\infty \lambda_n = \infty; \quad\text{and}\quad
 \sum_{n=1}^\infty | \; \lambda_{n+1} - \lambda_n \; | < \infty . $$
Then
 $\{ x_n \}$ converges strongly to a common fixed point of
 $\{ T(t) : t \geq 0 \}$.
\end{thm}


\begin{thebibliography}{99}

\bibitem{REF:Atsushiba_Takahashi2002_1}
 S. Atsushiba and W. Takahashi,
 {\it ``Strong convergence theorems for one-parameter nonexpansive
 semigroups with compact domains''},
 in Fixed Point Theory and Applications, Volume 3
 (Y. J. Cho, J. K. Kim and S. M. Kang Eds.),
 pp. 15--31,
 Nova Science Publishers, New York, 2002.

\bibitem{REF:Baillon1975}
 J. B. Baillon,
 {\it ``Un th\'eor\`eme de type ergodique pour les contractions non lin\'eaires
 dans un espace de Hilbert''},
 C. R. Acad.\ Sci.\ Paris, S\'er. A-B, {\bf 280} (1975), 1511--1514.

\bibitem{REF:Baillon1976}
 J. B. Baillon,
 {\it ``Quelques properi\'et\`es de convergence asymptotique pour les
 semigroupes de contractions impa{\`\i}res''},
 C. R. Acad.\ Sci.\ Paris, {\bf 283} (1976), 75--78.

\bibitem{REF:Banach1922_FundMath}
 S. Banach,
 {\it ``Sur les op\'erations dans les ensembles abstraits et leur application
 aux \'equations int\'egrales''},
 Fund.\ Math., {\bf 3} (1922), 133--181.

\bibitem{REF:Browder1965_ProcNAS_3}
 F. E. Browder,
 {\it ``Nonexpansive nonlinear operators in a Banach space''},
 Proc.\ Nat.\ Acad.\ Sci.\ USA, {\bf 54} (1965), 1041--1044.

\bibitem{REF:Browder1967_ARMA}
 F. E. Browder,
 {\it ``Convergence of approximates to fixed points of nonexpansive nonlinear
 mappings in Banach spaces''},
 Arch.\ Ration.\ Mech.\ Anal., {\bf 24} (1967), 82--90.

\bibitem{REF:Bruck1973_TransAMS}
 R. E. Bruck,
 {\it ``Properties of fixed-point sets of nonexpansive mappings in Banach
 spaces''},
 Trans.\ Amer.\ Math.\ Soc., {\bf 179} (1973), 251--262.

\bibitem{REF:Bruck1974_Pacific}
 R. E. Bruck,
 {\it ``A common fixed point theorem for a commuting family of nonexpansive
 mappings''},
 Pacific J.\ Math., {\bf 53} (1974), 59--71.

\bibitem{REF:Gohde1965}
 D. G\"ohde:
 {\it ``Zum Prinzip def kontraktiven Abbildung''},
 Math.\ Nachr., {\bf 30} (1965), 251--258.

\bibitem{REF:Halpern1967_BullAMS}
 B. Halpern,
 {\it ``Fixed points of nonexpanding maps''},
 Bull.\ Amer.\ Math.\ Soc., {\bf 73} (1967), 957--961.

\bibitem{REF:Hirano1982_JMSJapan}
 N. Hirano,
 {\it ``Nonlinear ergodic theorems and weak convergence theorems''},
 J.\ Math.\ Soc.\ Japan, {\bf 34} (1982), 35--46.

\bibitem{REF:Ishikawa1979_Pacific}
 S. Ishikawa,
 {\it ``Common fixed points and iteration of commuting nonexpansive
 mappings''},
 Pacific J.\ Math., {\bf 80} (1979), 493--501.

\bibitem{REF:Kirk1965_AMMonth}
 W. A. Kirk,
 {\it ``A fixed point theorem for mappings which do not increase distances''},
 Amer.\ Math.\ Monthly, {\bf 72} (1965), 1004--1006.

\bibitem{REF:KrasnoselskiiMA1955_Uspekhi}
 M. A. Krasnoselskii,
 {\it ``Two remarks on the method of successive approximations''}
 (in Russian),
 Uspehi Mat.\ Nauk {\bf 10} (1955), 123--127.

\bibitem{REF:Mann1953_ProcAMS}
 W. R. Mann,
 {\it ``Mean value methods in iteration''},
 Proc.\ Amer.\ Math.\ Soc., {\bf 4} (1953), 506--510.

\bibitem{REF:Miyadera_Kobayasi1982_NATMA}
 I. Miyadera and K. Kobayasi,
 {\it ``On the asymptotic behaviour of almost-orbits of nonlinear contraction
 semigroups in Banach spaces''},
 Nonlinear Anal., {\bf 6} (1982), 349--365.

\bibitem{REF:Reich1979_JMAA}
 S. Reich,
 {\it ``Weak convergence theorems for nonexpansive mappings''},
 J.\ Math.\ Anal.\ Appl., {\bf 67} (1979), 274--276.

\bibitem{REF:Shioji_Takahashi1998_NATMA}
 N. Shioji and W. Takahashi,
 {\it ``Strong convergence theorems for asymptotically nonexpansive semigroups
 in Hilbert spaces''},
 Nonlinear Anal., {\bf 34} (1998), 87--99.

\bibitem{REF:TS2002_JNCA}
 T. Suzuki,
 {\it ``Strong convergence theorem to common fixed points of two nonexpansive
 mappings in general Banach spaces''},
 J.\ Nonlinear Convex Anal., {\bf 3} (2002), 381--391.

\bibitem{REF:TS2003_ProcAMS}
 T. Suzuki,
 {\it ``On strong convergence to common fixed points of nonexpansive semigroups
 in Hilbert spaces''},
 Proc.\ Amer.\ Math.\ Soc., {\bf 131} (2003), 2133--2136.

\bibitem{REF:TSP_mann2uo_1_04}
 T. Suzuki,
 {\it ``Common fixed points of two nonexpansive mappings in Banach spaces''},
 to appear in Bull.\ Austral.\ Math.\ Soc.

\bibitem{REF:Suzuki-Takahashi}
 T. Suzuki and W. Takahashi
 {\it ``Strong convergence theorems of Mann's type for
 one-parameter nonexpansive semigroups
 in general Banach spaces''}, submitted.

\bibitem{REF:Takahashi_ybook}
 W. Takahashi,
 {\it ``Nonlinear Functional Analysis''},
 Yokohama Publishers,
 Yokohama, 2000.

\bibitem{REF:Wittmann1992_ArchMath}
 R. Wittmann,
 {\it ``Approximation of fixed points of nonexpansive mappings''},
 Arch.\ Math.\ (Basel), {\bf 58} (1992), 486--491.

\end{thebibliography}
\end{document}